\documentclass[12pt]{amsart}
\usepackage[mathscr]{eucal}
\newtheorem{thm}{Theorem}
\newtheorem{prop}[thm]{Proposition}
\newtheorem{cor}[thm]{Corollary}
\newtheorem*{thm*}{Theorem}
\newtheorem{lemma}[thm]{Lemma}

\begin{document}
\newcommand{\R}{{\mathbb R}}
\newcommand{\C}{{\mathbb C}}
\newcommand{\Z}{{\mathbb Z}}
\newcommand{\z}{{\mathfrak{Z}}}
\newcommand{\g}{{\mathfrak g}}
\newcommand{\Id}{{\mathbb I}}
\newcommand{\Ker}{\mathsf{Ker}}
\newcommand{\Out}{\mathsf{Out}}
\newcommand{\htg}{\hat{t}_{\gamma}}
\newcommand{\Aut}{\mathsf{Aut}}
\newcommand{\Ad}{\mathsf{Ad}}
\newcommand{\Inn}{\mathsf{Inn}}
\newcommand{\tr}{\mathsf{tr}}
\newcommand{\SU}{{\mathsf{SU}}(n)}
\newcommand{\Hom}{\mathsf{Hom}}
\newcommand{\hg}{\Hom(\pi ,G)/G}
\newcommand{\hpg}{\Hom(\pi ,G)}
\newcommand{\hz}{\Hom(\pi,\z)}
\newcommand{\Ou}{\Out(\pi)}
\newcommand{\Hh}{\mathcal{H}}
\newcommand{\W}{\mathcal{W}}

\title[Reducibility of Mapping Class Group Actions]
{The mapping class group acts reducibly on $\SU$-character varieties}
\author{William M.~Goldman}
\address{ Mathematics Department,
University of Maryland, College Park, MD  20742 USA  }
\email{ wmg@math.umd.edu }
\date{\today}
\begin{abstract}
When $G$ is a connected compact Lie group, and $\pi$ is a closed surface
group, then $\hg$ contains an open dense $\Ou$-invariant subset which
is a smooth symplectic manifold. This symplectic structure is $\Ou$-invariant
and therefore defines an invariant measure $\mu$, which has finite volume.
The corresponding unitary representation of $\Ou$ on $L^2(\hg,\mu)$ contains
no finite-dimensional subrepresentations besides the constants. 
This note gives a short proof that when $G=\SU$, the representation
$L^2(\hg,\mu)$ contains many other invariant subspaces.
\end{abstract}
\thanks{Partially supported by NSF grants DMS0405605 and DMS-0103889}
\maketitle
Let $G=\SU$ and $\pi$ be the fundamental group of a closed oriented surface 
$\Sigma$.
Let $\hpg$ be the space of  representations
$\pi\longrightarrow G$. The group $\Aut(\pi)\times\Aut(G)$ acts
on $\hpg$. Let $\hg$ be the quotient of $\hpg$ by $\{1\}\times\Inn(G)$.
Then $\Out(\pi) :=\Aut(\pi)/\Inn(\pi)$ acts on $\hg$.
The $\Ou$-action preserves a symplectic structure on $\hg$ which
determines a finite invariant smooth measure $\mu$ (on an invariant dense
open subset which is a smooth manifold). 
When $G$ is compact, 
the total measure is finite (Jeffrey-Weitsman~\cite{JW1,JW2}, 
Huebschmann~\cite{Huebschmann}).
There results a unitary representation of $\Ou$ on the Hilbert space 
\begin{equation*}
\Hh := L^2(\hg,\mu). 
\end{equation*}
Let $C\subset\Hh$ denote the subspace corresponding
to the constant functions. 

The following theorem is proved in Goldman~\cite{erg} for $n=2$ and
Pickrell-Xia~\cite{PickrellXia1,PickrellXia2} in general:

\begin{thm*} 
The only finite dimensional $\Ou$-invariant subspace of $\Hh$ is $C$.
\end{thm*}

According to \cite{erg}, the only finite-dimensional
invariant subspace of $\Hh$ consists of constants. Let $\Hh_0$ denote
the orthocomplement of the constants in $\Hh$, that is, the set of
$f\in\Hh$ such that $\int f d\mu = 0$. 

The goal of this note is an elementary proof of the following:

\begin{thm*} 
The representation of $\Ou$ on $\Hh_0$ is reducible.
\end{thm*}

In general for a compact connected Lie group $G$, the components of
$\hpg$ and $\hg$ are indexed by $\pi_1(G')$ where $G'$ is the
commutator subgroup. In that case the ergodicity/mixing results of
\cite{erg,PickrellXia1,PickrellXia2} 
imply that the only invariant finite dimensional
subspaces are subspaces of the space of locally constant functions, a vector
space of dimension $|vert\pi_1(G')\vert$. The method of proving reducibility
requires a nontrivial element of the center of $G$. 
For simplicity, we only discuss the case $G=\SU$ in this paper.

I am grateful to Steve Zelditch, Charlie Frohman, and Jeff Hakim
for helpful conversations, and to David Fisher for carefully reading
the manuscript and pointing out several corrections.

\section{The center of $\SU$}

Let $\z\cong \Z/n$ denote the center of $G$, the group consisting of all scalar
matrices $\zeta\Id$ where $\zeta^n=1$. 
Then $\hz \cong H^1(M;\Z/n)$ acts on $\hg$ by pointwise multiplication:
If $\rho\in\hpg$ and $u\in\hz$, then define the action $u\cdot\rho$ of $u$ on
$\rho$ by:
\begin{equation}\label{eq:defnaction}
\gamma\stackrel{u\cdot\rho}\longmapsto \rho(\gamma)u(\gamma). 
\end{equation}

Recall the definition \cite{Goldman1} of the symplectic structure on $\hg$.
Suppose that $\rho\in\hpg$ is an irreducible representation. By
Weil~\cite{Weil} (compare also Raghunathan~\cite{Raghunathan}), 
the Zariski tangent space to $\hpg$ at $\rho$
identifies with the space $Z^1(\pi,\g_{\Ad\rho})$ of cocycles where
$\g_{\Ad\rho}$ is the $\pi$-module defined by the composition
\begin{equation*}
\pi \overset{\rho}\longrightarrow G \overset{\Ad}\longrightarrow \Aut (\g) 
\end{equation*}
and $\g$ is the Lie algebra of $G$. 
Since $\rho$ is irreducible, $G$ acts locally freely and $\hg$ has a smooth
structure in a neighborhood of $[\rho]$. The tangent space
to the orbit $G\cdot\rho$ equals the space of coboundaries
$B^1(\pi,\g_{\Ad\rho})$. The tangent space to $\hg$ at $[\rho]$
identifies with the cohomology group
$H^1(\pi,\g_{\Ad\rho})$. 

Let $u\in\hz$. Since $\Ad(\z)$ is trivial, the action of $u$ induces
an identification of $\pi$-modules 
$\g_{\Ad\rho}\longrightarrow\g_{\Ad(u\cdot\rho)} $
and hence of tangent spaces

\begin{align*}
T_{[\rho]} \hg & =  
H^1(\pi,\g_{\Ad\rho}) \\ & \quad \longrightarrow \quad
T_{[u\cdot\rho]} \hg  = H^1(\pi,\g_{\Ad(u\cdot\rho)}).
\end{align*}

The symplectic form $\omega_{\rho}$ at $[\rho]$ is defined by the cup
product 
\begin{equation*}
H^1(\pi,\g_{\Ad\rho}) \times H^1(\pi,\g_{\Ad\rho})
\longrightarrow H^2(\pi;\R)
\end{equation*}
using the pairing of $\pi$-modules induced by 
an $\Ad$-invariant nondegenerate symmetric bilinear form on $\g$
as a coeefficient pairing. Evidently the symplectic form $\omega$,
and therefore the corresponding measure $\mu$, are $\hz$-invariant.

\section{Trace functions}

Suppose that $\gamma\in\pi$. The function
\begin{align*}
\hpg& \longrightarrow \C \\ 
\rho & \longmapsto \tr(\rho(\gamma))
\end{align*}
is $\Inn(G)$-invariant and defines a function
\begin{equation*}
\hg\xrightarrow{t_{\gamma}}\C. 
\end{equation*}

We extend this definition to products of more than one element of $\pi$:
Let $\gamma = (\gamma_1,\dots,\gamma_s)\in \pi^s$.
Define the {\em trace function\/} of $\gamma$ as the product:
\begin{equation*}
t_\gamma ([\rho]) := t_{\gamma_1}([\rho]) \dots  t_{\gamma_s}([\rho]) 
\end{equation*}

Using the definition \eqref{eq:defnaction} of the action of $\hz$ on $\hg$,
the action of $u\in\hz$ on the trace function $t_\gamma$ is given by
$u\cdot t_\gamma$, defined by:
\begin{equation}\label{eq:funact}
\rho \stackrel{u\cdot t_\gamma}\longmapsto 
u(\gamma)^{-1}t_\gamma \big(\rho\big)  
\end{equation}
since
\begin{align*}
\big(u\cdot t_\gamma\big)\big( \rho \big) & :=
t_\gamma(u^{-1}\cdot \rho) 
\\ & = \tr\big( (u^{-1}\cdot \rho)(\gamma)\big) 
\\ & = \tr\big( u(\gamma)^{-1}\rho(\gamma) \big)
\\ & = u(\gamma)^{-1} \tr(\rho(\gamma))  
\\ & = u(\gamma)^{-1}t_\gamma \big(\rho\big)  
\end{align*}

Since $\z \cong \Z/n$, the evaluation $u(\gamma)$ depends only
on the {\em total homology class\/}
\begin{equation*}
[\gamma] :=  [\gamma_1] + \dots + [\gamma_s] \in H_1(\Sigma;\Z/n),
\end{equation*}
defined as the sum of the $\Z/n$-homology classes of the $\gamma_i$.
The evaluation of $u\in\hz$ on $\gamma$ is just the natural pairing
\begin{align*}
\hz \times H_1(\Sigma;\Z/n) & \longrightarrow \Z/n \cong Z \\
 (u, [\gamma]) & \longmapsto u\cdot[\gamma]
\end{align*}
where $[\gamma]$ is the total homology class of $\gamma$.

%

\begin{lemma} \label{lemma}
If $[\gamma]\neq 0 \in H_1(\Sigma;\Z/n)$, then
\begin{equation*}
\int t_\gamma d\mu = 0.
\end{equation*}
\end{lemma}
\begin{proof}
Since the homology class of $\gamma$ in $H_1(\Sigma,\Z/n)$ is nonzero,
$u(\gamma)\neq 1$ for some $u\in\hz$. Since
\big($u^{-1}\big)_*\mu = \mu$, 
\begin{equation*}
\int t_\gamma \,d\mu \;=\; \int  (t_\gamma\circ u^{-1})
d\big(u^{-1}\big)_*\mu  =
u(\gamma)^{-1}\int t_\gamma d\mu
\end{equation*}
by  \eqref{eq:funact}. Since $u(\gamma)^{-1}\neq 1$, the integral is zero
as claimed.
\end{proof}

An element $\gamma\in\pi$ is {\em nonseparating simple\/} if
it contains a nonseparating simple loop which is homotopically nontrivial
and not homotopic to $\partial \Sigma$.
Since such an element has nontrivial homology class
$[\gamma]\in H_1(\Sigma,\Z/n)$,  Lemma~\ref{lemma} implies:

\begin{cor}\label{cor}
 Suppose $\gamma$ is nonseparating simple. Then $t_\gamma\in\Hh_0$.
\end{cor}

For $n = 2$, Frohman and Kania-Bartoszynska~\cite{Frohman} have calculated
$\int t_\gamma d\mu$ for separating simple loops $\gamma$ in terms
of Bernoulli numbers and $6j$-symbols. 

\section{Invariant subspaces}

Let $\W\subset \Hh_0$ denote the closure
of the span of all $t_\alpha$, where $\alpha$ is nonseparating simple loops.
(Alternatively, let $\alpha$ run over all elements of $\pi$ with nonzero
homology class modulo $n$.) Clearly $\W$ is $\Ou$-invariant.

\begin{prop}
$\W$ is a proper subspace of $\Hh _0$.
\end{prop}
\begin{proof}
Choose $\gamma\in\pi^s$ such that:
\begin{itemize}
\item $t_\gamma$ is not constant;
\item $\gamma$ has trivial homology class in $H_1(\Sigma,\Z/n)$.
\end{itemize}
Such elements are easy to find: for example if $\gamma$ consists of one 
nontrivial separating simple loop, or if $\gamma$ consists of one
$n$-th power, or if $\gamma=(\gamma_1,\dots,\gamma_1)$ where $s$ is a positive
multiple of $n$ and $\gamma_1$ is nontrivial.

Normalize $t_\gamma$ so that it lies in $\Hh_0$:
\begin{equation*}
\htg :=  t_\gamma - \int t_\gamma d\mu.
\end{equation*}

We claim that $\htg$ is a nonzero vector in $\Hh_0$ orthogonal
to $\W$. It is nontrivial since $t_\gamma$ is nonconstant. 
Furthermore if $\alpha\in\pi$ has nonzero homology class 
$[\alpha]\in H_1(\Sigma,\Z/n)$ (for example if it nonseparating simple), then
$t_\alpha$ is orthogonal to $\htg$. 

Since $\rho(\alpha)$ is unitary, $\overline{t_\alpha} = t_{\alpha^{-1}}$.
Thus
\begin{align*}
\langle \htg ,t_\alpha\rangle & = 
\int \htg \overline{t_\alpha} d\mu   =
\int \htg t_{\alpha^{-1}} d\mu  \\ & =
\int t_\gamma t_{\alpha^{-1}} d\mu  \text{\quad~(by Corollary~\ref{cor})}\\ & =
 \int t_\eta  d\mu  
\end{align*}
where $\eta = (\gamma_1,\dots,\gamma_s,\alpha^{-1})$. 
Since 
\begin{equation*}
[\eta] = [\gamma] - [\alpha] \neq 0,  
\end{equation*}
Lemma~\ref{lemma} implies that $\langle \htg,t_\alpha\rangle = 0$, as
desired.
\end{proof}

When $\pi$ is a free group and $n=2$, then ergodicity and weak-mixing
of the action of $\Ou$ on $\hg$ is proved in Goldman~\cite{outfn}.
Exactly the same argument as above (with $\W$ replaced by the closure
of the span of $t_\gamma$ where $\gamma$ is an element of a free generating
set of $\pi$) shows that the unitary representation of $\Ou$ on 
$L^2(\hg,\mu)$ contains invariant (necessarily infinite-dimensional)
subspaces other than the constants.

%
%
%
%

%

\makeatletter \renewcommand{\@biblabel}[1]{\hfill#1.}\makeatother
 \renewcommand{\bysame}{\leavevmode\hbox to3em{\hrulefill}\,}


\begin{thebibliography}{10}

\bibitem{Frohman} 
Frohman, C. and Kania-Bartoszynska, J.,
  {\em Shadow world evaluation of the Yang-Mills measure,\/}
  Algebraic and Geometric Topology, {\bf 4 \/} (2001), 311--332.


\bibitem{Goldman1} 
Goldman, W.,	{\em The symplectic nature of fundamental groups of 
		surfaces,\/} 
		Adv.\ Math.\ {\bf 54} (1984), 200--225.



\bibitem{erg}	
		\bysame,
		{\em Ergodic theory on moduli spaces,\/}
		Ann.\ Math.\  {\bf 146} (1997), 475--507.

\bibitem{outfn} \bysame,
            {\em An ergodic action of the outer automorphism group 
	      of a free group,\/} {\tt math.DG/0506401} (submitted)


\bibitem{Huebschmann}  Huebschmann, J.,
  {\em Symplectic and Poisson structures of certain moduli spaces, \/}
  Duke Math J. {\bf 80} (1995) 737-756.
                                                                            
\bibitem{JW1} Jeffrey,L. and  Weitsman, J.,
{\em  Bohr-Sommerfeld orbits and the Verlinde dimension formula, \/}
Commun. Math. Phys. 150 (1992) 593-630
                                                                              
\bibitem{JW2} \bysame,
{\em Toric structures on the moduli space of flat connections on a Riemann
surface: Volumes and the moment map,\/}
Adv.\ Math. {\bf 109,} 151-168 (1994).

\bibitem{PickrellXia1} Pickrell, D.\ and Xia, E.,
		{\em Ergodicity of Mapping Class Group Actions on 
		Representation Varieties, I. Closed Surfaces,\/}
		Comment.\ Math.\ Helv.\ {\bf 77} (2001), 339--362.


\bibitem{PickrellXia2} \bysame,
		{\em Ergodicity of Mapping Class Group Actions on 
		Representation Varieties, II. Surfaces with Boundary,\/}
 		Transformation Groups {\bf 8} (2003), no.\ 4, 397--402.


\bibitem{Raghunathan}	
		Raghunathan, M.,
		``Discrete Subgroups of Lie Groups,''
		Ergebni\ss e der Math. {\bf 58},
		Springer-Verlag Berlin-Heidelberg-New York (1972).
		
\bibitem{Weil}	Weil, A.,
		{\em Remarks on the cohomology of groups,\/}
		Ann.\ Math.\ {\bf 80} (1964), 149--157.

\end{thebibliography}
\end{document}